\documentclass{article}
\usepackage{graphicx}
\usepackage{amsmath}
\usepackage{amssymb}
\usepackage{amsthm}
\usepackage{mathrsfs}
\usepackage{enumerate}

\newtheorem{thm}{Theorem}

\newtheorem{cor}{Corollary}

\newtheorem{example}{Example}

\title{A Generalization of ``Existence and Behavior of the Radial Limits of a Bounded Capillary
       Surface at a Corner''}
\author{Julie N. Crenshaw, Alexandra K. Echart \& Kirk E. Lancaster
                                  \\
                       Department of Mathematics, Statistics \& Physics \\
                            Wichita State University \\
                            Wichita, Kansas, 67260-0033}

\def\Real{{\rm I\hspace{-0.2em}R}}
\def\Natural{{\rm I\hspace{-0.2em}N}}

\DeclareMathOperator*{\esslimsup}{ess\,lim\,sup}
\DeclareMathOperator*{\essliminf}{ess\,lim\,inf}

\begin{document}
\maketitle

\vspace{5mm}

\begin{abstract}
The principle existence theorem (i.e. Theorem 1) of ``Existence and Behavior of the Radial Limits of a Bounded 
Capillary Surface at a Corner'' (Pacific J.~Math. {\bf Vol. 176}, No. 1 (1996), 165-194) is extended to the case 
of a contact angle $\gamma$  which is  not bounded away from $0$  and $\pi$  (and depends on position in a bounded 
domain $\Omega\in \Real^{2}$  with a convex corner at ${\cal O}=(0,0)$).  The lower bound on the size of 
``side fans'' (i.e. Theorem 2 in the above paper) is extended to case of such contact angles for convex and 
nonconvex corners.  
\end{abstract}

\section{Introduction and Theorems}

Consider the capillary problem
\begin{eqnarray}
\label{eq:abasic}
Nf & = & \kappa f + \lambda  \mbox{  \ in \ } \Omega    \\
Tf \cdot {\bf \nu} & = & \cos \gamma  \mbox{ \ on \ } \partial \Omega
 \label{eq:bbasic}
\end{eqnarray}
where $\Omega$ is a region in ${\Real}^{2}$ with a corner at ${\cal O}$, ${\cal O} \in \partial\Omega, \: Nf = \nabla \cdot Tf, \: 
Tf = \frac{\nabla f}{\sqrt{1 + |\nabla  f|^{2}}},$  
$\kappa$ and $\lambda$ are constants, ${\bf \nu}$ is the exterior unit normal on
$\partial \Omega,$ and $\gamma = \gamma (s)$ is a function of position on
$\partial \Omega, 0 \leq \gamma (s) \leq \pi.$  The surface $z = f(x,y)$
describes the shape of the static liquid-gas interface in a vertical cylindrical
tube of cross-section $\Omega$; see  \cite{Finn:86,LS1} for background.  
 
We are interested in the behavior of solutions to (\ref{eq:abasic}), (\ref{eq:bbasic}) in a 
neighborhood of a corner point of the boundary.  We take the corner point to be ${\cal O}=(0,0).$  
Let $\Omega^{*} = \Omega \cap B_{\delta^{*}}({\cal O})$, where  $B_{\delta^{*}}({\cal O})$ is the ball of radius 
$\delta^{*}$ about ${\cal O}$.  Polar coordinates relative to ${\cal O}$ will be denoted by $r$ and $\theta$.
We assume that $\partial \Omega$ is piecewise smooth and that $\partial \Omega \cap B_{\delta^{*}}({\cal O})$ consists of two arcs
${\partial}^{+}\Omega^{*}$  and $\partial^{-}\Omega^{*}$, whose tangent lines approach the lines $L^{+}: \:  \theta = \alpha$
and $L^{-}: \: \theta = - \alpha$, respectively, as the point ${\cal O}$ is approached.
The points where $\partial B_{\delta^{*}}({\cal O})$ intersect $\partial \Omega$  are labeled $A$ and $B$; 
also, $\Gamma^{*} = \partial B_{\delta^{*}}({\cal O})\cap \overline{\Omega}.$
Set 
\[
\Omega_{\infty} = \{\left(r\cos(\theta),r\sin(\theta)\right) \ : \ r>0, -\alpha<\theta<\alpha\}.
\]
Let $(x^{+}(s),y^{+}(s))$  be an arclength parametrization of ${\partial}^{+}\Omega^{*}$  and $(x^{-}(s),y^{-}(s))$  
be an arclength parametrization of ${\partial}^{-}\Omega^{*},$  each measured from the corner at ${\cal O},$ 
so that $\left(x^{\pm}(0),y^{\pm}(0)\right)=(0,0).$  
Let $(x^{+}_{*}(s),y^{+}_{*}(s))$  be an arclength parametrization of 
$\partial^{+}\Omega_{\infty} = \{\left(r\cos(\alpha),r\sin(\alpha)\right) \ : \ r\ge 0\}$   
and  $(x^{-}_{*}(s),y^{-}_{*}(s))$  be an arclength parametrization of 
$\partial^{-}\Omega_{\infty} = \{\left(r\cos(-\alpha),r\sin(-\alpha)\right) \ : \ r\ge 0\},$
each measured from the corner at $(0,0).$ 
Define 
\[
\gamma^{+}(s)=\gamma\left(x^{+}(s),y^{+}(s)\right)\ \ \ \ {\rm and} \ \ \ \ \gamma^{-}(s)=\gamma\left(x^{-}(s),y^{-}(s)\right).
\]
For $0 \leq \alpha \le \pi/2$, the corner will be said to be {\em convex} and for $\pi/2 < \alpha \leq \pi$, 
the corner will be said to be {\em nonconvex}.

In \cite{LS1}, the existence of radial limits of a bounded solution $f$ to (\ref{eq:abasic}) that satisfies 
(\ref{eq:bbasic}) on the smooth portions of $\partial \Omega$  is proven provided that $\gamma$ was bounded away from $0$ and $\pi$, 
and for a convex corner an additional condition is satisfied coupling $\gamma^{+}$  and $\gamma^{-}$.  
In this paper, we eliminate the requirement that $\gamma$  is bounded away from $0$ and $\pi;$  
an additional condition must still be satisfied at a convex corner.  
The radial limits of $f$ will be denoted by 
\[
Rf(\theta)=\lim_{r \rightarrow 0^{+}} f(r\cos \theta, r\sin \theta), -\alpha < \theta < \alpha
\]
and $Rf(\pm \alpha)=\lim_{\partial^{\pm}\Omega^{*}\ni {\bf x} \rightarrow {\cal O}} f({\bf x}), {\bf x} = (x,y)$, 
which are the limits of the boundary values of $f$ on the two sides of the corner if these exist.

\begin{thm}
\label{ONE}
Let $f$ be a bounded solution to (\ref{eq:abasic}) satisfying (\ref{eq:bbasic}) on $\partial^{\pm}\Omega^{*} \setminus \{{\cal O}\}$ 
be discontinuous at ${\cal O}.$   
\begin{itemize}
\item[(a)]  If $\alpha > \pi/2$ then $Rf(\theta)$ exists for all $\theta \in (-\alpha,\alpha).$   
\item[(b)]If $\alpha \le \pi/2$ and there exist constants $\underline{\gamma}^{\, \pm},
\overline{\gamma}^{\, \pm}, 0 \le \underline{\gamma}^{\, \pm} \leq \overline{\gamma}^{\, \pm} \le \pi,$   satisfying
\[ 
\pi - 2\alpha  < \underline{\gamma}^{+} + \underline{\gamma}^{-} \le  
\overline{\gamma}^{\, +} +  \overline{\gamma}^{\, -} < \; \pi + 2\alpha 
\]
such that $\underline{\gamma}^{\pm}\leq \gamma^{\pm}(s) \leq \overline{\gamma}^{\,  \pm}$
for all $s\in (0,s_{0}),$ for some $s_{0}>0$, then $Rf(\theta)$ exists for
all $\theta \in (-\alpha, \alpha)$.  
\end{itemize}
Furthermore, in either case, $Rf(\theta)$  is a continuous
function on $(-\alpha, \alpha)$ which behaves in one of the following ways:
 
\noindent {\bf (i)} $Rf(\theta)$  is a constant function of $\theta$  and $f$  has a nontangential limit at ${\cal O}.$  

\noindent {\bf (ii)} There exist $\alpha_{1}$ and $\alpha_{2}$ so that $-\alpha \leq \alpha_{1}
< \alpha_{2} \leq \alpha$ and $Rf$ is constant on $(-\alpha, \alpha_{1}]$ and
$[ \alpha_{2}, \alpha)$ and strictly increasing or strictly decreasing on
$[ \alpha_{1}, \alpha_{2}]\cap (-\alpha,\alpha)$.  Label these case(I) and case(D), respectively.
 
\noindent {\bf (iii)} There exist $\alpha_{1}, \alpha_{L}, \alpha_{R}, \alpha_{2}$ so that
$-\alpha \leq \alpha_{1} < \alpha_{L} < \alpha_{R} < \alpha_{2} \leq \alpha,
\alpha_{R}= \alpha_{L} + \pi$, and $Rf$ is constant on $(-\alpha, \alpha_{1}],
[ \alpha_{L}, \alpha_{R}]$, and $[ \alpha_{2}, \alpha)$ and either increasing
on $[ \alpha_{1}, \alpha_{L}]\cap (-\alpha,\alpha)$ and decreasing on $[ \alpha_{R}, \alpha_{2}]\cap (-\alpha,\alpha)$ or
decreasing on $[ \alpha_{1}, \alpha_{L}]\cap (-\alpha,\alpha)$ and increasing on $[ \alpha_{R},
\alpha_{2}]\cap (-\alpha,\alpha)$.  Label these case (ID) and case (DI), respectively.
\end{thm}
 
\noindent 
In Theorem 1 of \cite{LS1} and Theorem \ref{ONE} above, the existence of two intervals $(-\alpha,\alpha_{1}]$  and 
$[\alpha_{2},\alpha)$  on which $Rf(\cdot)$  is constant (i.e. ``side fans'') is established but the relationship between 
the sizes of these side fans and the contact angle is unclear.  Theorem 2 of \cite{LS1} establishes lower bounds on these sizes 
when the 
\[
\lim_{\partial^{+}\Omega\ni (x,y)\to {\cal O}} \gamma(x,y) = \gamma_{0}^{+} \ \ \ {\rm and} \ \ \ 
\lim_{\partial^{-}\Omega\ni (x,y)\to {\cal O}} \gamma(x,y) = \gamma_{0}^{-}
\]
are assumed to exist.  (In \cite{CFC,NCFC}, these lower bounds were shown to be the actual sizes of the side fans.) 
What happens if the limits of $\gamma$  at ${\cal O}$  do not exist?  Theorem \ref{TWO} and Corollary \ref{COR1} 
provide lower bounds in this situation. 
 
For $0<b<1,$  define 
\[
A^{\pm}_{I}(b) = \liminf_{\epsilon\downarrow 0} \frac{1}{\epsilon} \int_{0}^{b\epsilon} \ \cos\left(\gamma^{\pm}(t)\right) dt
\ \ \ {\rm and} \ \ \ 
A^{\pm}_{S}(b) = \limsup_{\epsilon\downarrow 0} \frac{1}{\epsilon} \int_{0}^{b\epsilon} \ \cos\left(\gamma^{\pm}(t)\right) dt. 
\]
Notice that $b \cos\left(\limsup_{t\downarrow 0} \gamma^{\pm}(t)\right) \le A^{\pm}_{I}(b)  \le A^{\pm}_{S}(b)
\le b \cos\left(\liminf_{t\downarrow 0} \gamma^{\pm}(t)\right).$  

\begin{thm}  
\label{TWO}
Let $f$ be a bounded solution to (\ref{eq:abasic}) satisfying (\ref{eq:bbasic}) on 
$\partial^{\pm}\Omega^{*} \setminus \{{\cal O}\}$  which is discontinuous at ${\cal O}.$ 
Assume  $Rf(\theta)$ exists for all $\theta \in (-\alpha,\alpha).$   Then: 
\begin{itemize}
\item [(a)]  $Rf(\theta)$  is a continuous function on $(-\alpha, \alpha)$ which behaves as described 
in {\bf (i)}, {\bf (ii)} or {\bf (iii)}  of Theorem \ref{ONE}.  
\item[(b)] There exist fans of constant radial limits adjacent to each tangent direction at ${\cal O}$  and lower bounds on the sizes of 
these side fans exist.
\end{itemize}
In terms of the cases labeled in Theorem 1, 
the sizes of the side fans $\beta^{-}=\alpha_{1} + \alpha$  and $\beta^{+}=\alpha - \alpha_{2}$  satisfy the 
following conditions:     
\begin{itemize}
\item[(1.)]  $A_{I}^{+}\left(\frac{\sin(\lambda-\beta^{+})}{\sin(\lambda)}\right) + \frac{\sin(\beta^{+})}{\sin(\lambda)}\ge 1$ 
for all $\lambda\in (\beta^{+},\pi)$  for (I) and (DI).
\item[(2.)]  $A_{I}^{-}\left(\frac{\sin(\lambda-\beta^{-})}{\sin(\lambda)}\right) + \frac{\sin(\beta^{-})}{\sin(\lambda)}\ge 1$ 
for all $\lambda\in (\beta^{-},\pi)$  for  (D) and (DI).
\item[(3.)] $1+A_{S}^{-}\left(\frac{\sin(\lambda-\beta^{-})}{\sin(\lambda)}\right)  \le \frac{\sin(\beta^{-})}{\sin(\lambda)}$ 
for all $\lambda\in (\beta^{-},\pi)$  for  (I) and (ID).
\item[(4.)] $1+A_{S}^{+}\left(\frac{\sin(\lambda-\beta^{+})}{\sin(\lambda)}\right)  \le \frac{\sin(\beta^{+})}{\sin(\lambda)}$ 
for all $\lambda\in (\beta^{+},\pi)$  for (D) and (ID).
\end{itemize}
\end{thm}

\section{Proofs of Theorems \ref{ONE} \& \ref{TWO}} 
The proof of Theorem \ref{ONE} follows that established in \cite{Lan1} and \cite{EL:86} in which (i) the graph of the solution in $\Omega\times\Real$  is 
represented in isothermal coordinates, (ii) comparison arguments are used to prove that the component functions of the isothermal parametrization 
of the graph are uniformly continuous and so extend to be continuous on the closure of the parameter domain, (iii) boundary regularity theory 
(e.g. \cite{Heinz}) is used to prove that radial limits exist for almost every direction, (iv) cusp solutions are excluded (e.g. \cite{EchartLancaster1}) 
and (v) the behavior of the radial limit function is determined.  
The only step which does not follow from previous work is (ii) and so the proof of Theorem \ref{ONE} comes down to establishing (ii).  
The proof of Theorem \ref{TWO} follows from standard ``blow up'' arguments.

\subsection{Proof of Theorem \ref{ONE}:}  
When $\alpha>\frac{\pi}{2},$  Theorem \ref{ONE} is a consequence of \cite{NoraKirk1}.  
Suppose now that $\alpha\le \frac{\pi}{2}.$ 
Since  $f$  is bounded and the prescribed mean curvature is $H(x,y,z)=\kappa z+\lambda,$  there exist $M_{1} \in (0, \infty)$ and 
$M_{2} \in [0, \infty)$ such that 
\begin{equation}
\label{Bounds}
\sup_{(x,y) \in \Omega} |f(x,y)| \leq M_{1} \mathrm{ \ \ and \ } \sup_{(x,y) \in \Omega} |H(x, y, f(x,y))| \leq M_{2}.
\end{equation}  
In \S 2.1 of \cite{NoraKirk2}, a specific torus is constructed which depends solely on $M_{2}$  and which is used as a comparison surface; 
one should compare this with, for example, \cite{LS1} where several types of comparison surfaces are used or \cite{NoraKirk1} where an unduloid 
is used as a comparison surface.  We shall use this torus as our comparison surface here. 
We will denote by $q$ the denote the modulus of continuity of the function $h^{-}$  whose graph is the set $\mathcal{T}$  which is 
the inner half of a torus with axis of symmetry $\{(2, y, 0): y \in\mathbb{R}^{2}\}$, major radius $R_{0} = 2$, and minor radius $r_{0};$  here  
\begin{equation}
\label{R0}
r_{0} = \begin{cases}
1 & \mathrm{if \ } M_{2} = 0 \\
\frac{1}{M_{2}} + 1 - \sqrt{\left(\frac{1}{M_{2}}\right)^{2} + 1} & \mathrm{if \ } M_{2} > 0.
\end{cases}
\end{equation}
Then $q$ is also the modulus of continuity of functions (i.e. $h^{+}, h_{\beta}^{-}$, $h_{\beta}^{+}$) whose graphs are obtained by rotations and 
translations in the horizontal plane of $\mathcal{T}$  (see \cite{NoraKirk2}, page 59). 

Let $\mathscr{S}_{0} = \mathrm{gra}(f) = \{(x,y,f(x,y)): (x,y) \in \Omega^{*}\}$ and allow $\mathscr{S}$ to be the closure of 
$\mathscr{S}_{0}$  in $\mathbb{R}^{3}$. 
As in \S 2.2 of \cite{NoraKirk2}, there exists an isothermal parametrization $Y: E \rightarrow \mathbb{R}^{3}$ given by 
\[
Y(u,v) = (a(u,v), b(u,v), c(u,v))
\]
such that $Y(\bar{E}) = \mathscr{S}$,  $Y(E) = \mathscr{S}_{0}$, and ($a_1$) - ($a_5$) of \cite{NoraKirk2} hold, where 
$E = B_{1}(\mathcal{O}) = \{(u,v): u^{2} + v^{2} < 1\}.$  
By $(a_{2})$  of \cite{NoraKirk2}, if we let $G(u,v) = (a(u,v), b(u,v))$  for $(u,v) \in E$, then $G \in C^{0}(\bar{E})$.  
From ($a_{3})$  of \cite{NoraKirk2}, there exists a connected arc $\sigma \subset \partial E$ that $Y$ maps strictly monotonically onto 
$\{(x,y, f(x,y)): (x,y)\in\partial\Omega^{*} \setminus \{\mathcal{O}\}\}.$   
Let the endpoints of $\sigma$ be denoted $\textbf{o}_{1}$ and $\textbf{o}_{2}.$  
There exists points $\textbf{a}, \textbf{b} \in \sigma$ such that $G(\textbf{a}) = A,$  $G(\textbf{b}) = B,$  
$G$ maps the arc $\textbf{o}_{2} \textbf{a}$ onto $\partial^{-}\Omega$  and $G$ maps the arc $\textbf{o}_{1}\textbf{b}$ onto 
$\partial^{+}\Omega.$   
We must consider the two cases:  
\begin{enumerate}[(A)]
\item $\textbf{o}_{1} = \textbf{o}_{2}$
\item $\textbf{o}_{1} \neq \textbf{o}_{2}$
\end{enumerate}
 
\noindent {\bf Assume first that (A) holds.}  Set $\textbf{o} = \textbf{o}_{1} = \textbf{o}_{2}$. 
We wish to prove that $c$ is uniformly continuous on $E$ and hence $c$ extends to be continuous on $\bar{E}.$ 
If so, then the existence and behavior of the radial limits of $f$  follows as in \cite{NoraKirk2,LS1}. 
There are three possible cases: \begin{enumerate}[(i)] 
\item $\underline{\gamma}^{-} > 0$ and $\bar{\gamma}^{-} < \pi$
\item $\underline{\gamma}^{+} > 0$ and $\bar{\gamma}^{+} < \pi$
\item $(\underline{\gamma}^{-} = 0$ or $\bar{\gamma}^{-} = \pi)$ and $(\underline{\gamma}^{+} = 0$ or $\bar{\gamma}^{+} = \pi)$
\end{enumerate}

\noindent {\bf Case (i):} Let $\lambda_{1} = \underline{\gamma}^{+}, \lambda_{2} = \bar{\gamma}^{+}, \gamma_{2} = \underline{\gamma}^{-}.$ 
We observe that $\lambda_{2} = \bar{\gamma}^{+} < \pi + 2\alpha - \bar{\gamma}^{-}$, 
$\lambda_{1} = \underline{\gamma}^{+} > \pi - 2\alpha - \underline{\gamma}^{-}$, and so $\lambda_{2} - \lambda_{1} < 4\alpha$. 
We wish to use the argument in the proof of Theorem 2 of \cite{NoraKirk2}.  
Since $\pi-2\alpha-\lambda_{1}<\gamma_{2}<\pi+2\alpha-\lambda_{2},$  we can choose $\tau_{1},\tau_{2}\in (0,\pi)$  
such that $\tau_{1}\in (\pi-2\alpha-\lambda_{1}, \gamma_{2})$   and $\tau_{2}\in (\gamma_{2},\pi+2\alpha-\lambda_{2}).$
Set $\beta_{1}=\frac{\pi}{2}-\tau_{1}$  and  $\beta_{2}=\tau_{2}-\frac{\pi}{2}.$  
With these choices of $\beta_{1}$  and $\beta_{2},$   notice that 
\[
T\left( h^{-}\circ T_{\beta_{1}}\right)(x_{1},0)\cdot (0,-1) = \cos\left(\tau_{1}\right)
> \cos(\gamma_{2}), \ \  {\rm for} \ 0<x_{1}< 2-r_{0}
\]
and 
\[
T\left( h^{+}\circ T_{\beta_{2}}\right)(x_{1},0)\cdot (0,-1) = \cos\left(\tau_{2}\right)<\cos(\gamma_{2}),\ \ {\rm for} \ 0<x_{1}< 2-r_{0}
\]
(see \cite{NoraKirk2}, page 59). 
This implies that for $\delta_{1}=\delta_{1}(\beta_{1},\beta_{2})>0$  small enough and 
${\bf x}\in\partial^{-}\Omega$  with $|{\bf x}|<\delta_{1},$  we have  
\begin{equation}
\label{ContactA}
T\left( h^{-}_{\beta_{1}}\right)({\bf x})\cdot \vec \nu({\bf x}) > \cos(\gamma({\bf x}))\ \ \ {\rm and} \ \ \ 
T\left( h^{+}_{\beta_{2}}\right)({\bf x})\cdot \vec \nu({\bf x}) < \cos(\gamma({\bf x})).
\end{equation}
Since $\beta_{1},\beta_{2}\neq \pm \frac{\pi}{2},$  there exists $R=R(\beta_{1},\beta_{2})>0$  such that 
$B_{R}({\cal O})\cap \Omega^{*}\subset \Delta_{\beta_{1}}\cap \Delta_{\beta_{2}},$  
where $\Delta_{\beta}$  is as in \S 2.1 of \cite{NoraKirk2}.
For each $\delta \in (0,1),$  allow 
\begin{equation}
p(\delta) = \sqrt{\frac{8\pi M_{0}}{\ln \frac{1}{\delta}}}
\end{equation} 
where $M_{0}$ is the area of $S_{0}.$   

Let $\epsilon > 0.$  Choose $\delta > 0$ such that 
\[
\sqrt{\delta} < \min\{ ||\textbf{o} - \textbf{a}||, ||\textbf{o} - \textbf{b}||\}, \ \ 
p(\delta) < \delta_{1}(\beta_{1},\beta_{2}), \ \ p(\delta) < R(\beta_{1}, \beta_{2}) \ \ 
{\rm and} \ p(\delta) + q (p(\delta)) < \frac{\epsilon}{2}. 
\]
Let ${\bf w_{0}}=(u_{0},v_{0})\in E.$  From the Courant-Lebesgue lemma, there exists a $\rho(\delta) \in (\delta, \sqrt{\delta})$  
such that the arclength $l_{\rho(\delta)}$ of $C_{\rho(\delta)}'$  is less than $p(\delta),$  where 
$C_{\delta} = \{{\bf w} \in E: ||{\bf w} - {\bf w_{0}}|| = \delta\}$  and $C_{\delta}' = Y(C_{\delta}).$  
Set $B_{\delta} = \{{\bf w} \in E: ||{\bf w} - {\bf w_{0}}|| < \delta\}$  and $B_{\delta}' = Y\left(B_{\delta}\right).$ 
Then, for $\textbf{w} \in C_{\rho(\delta)}',$ there exist functions 
\begin{equation}
b^{+} (x,y) = f(\textbf{w}) + p(\delta) + h_{\beta_{1}}^{-} (x,y) \mathrm{ \ for \ } (x,y) \in \Delta_{\beta_{1}}
\end{equation}
\begin{equation}
b^{-} (x,y) = f(\textbf{w}) - p(\delta) - h_{\beta_{2}}^{+} (x,y) \mathrm{ \ for \ } (x,y) \in \Delta_{\beta_{2}}
\end{equation} 
where $\beta_{1} = \frac{\pi}{2} - \tau_{1}$ and $\beta_{2} = \tau_{2} - \frac{\pi}{2}$. 
From (10) of \cite{NoraKirk2},  we have that $\mathrm{div}(b^{+}) \leq - M_{2}$ in $\Delta_{\beta_{1}}$ and 
$\mathrm{div}(b^{-}) \geq M_{2}$ in $\Delta_{\beta_{2}}.$  
So in $\Omega \cap \Delta_{\beta_{1}}$, $\mathrm{div} (Tb^{+}) \leq \mathrm{div} (Tf)$. 
On $\partial^{-}\Omega \cap B_{\delta_{1}}(\mathcal{O})$, $Tb^{+} \cdot \nu \geq Tf \cdot \nu$. 
As in the proof of Theorem 2 of \cite{NoraKirk2}, 
\begin{equation}
\label{R4}
f(x,y) < b^{+}(x,y) \ \ \ \ {\rm for} \ \ (x,y) \in \Delta_{\beta_{1}} \cap B_{\rho(\delta)}',
\end{equation} 
where $B_{\rho(\delta)}' = Y\left(B_{\rho(\delta)}\right).$   This follows since $Tb^{+} \cdot \nu \geq Tf \cdot \nu$  
on $\partial^{+}\Omega \cap B_{\delta_{2}}(\mathcal{O})$  by (15) of \cite{NoraKirk2} if $\tau_{1}+2\alpha\le \pi$
and no boundary condition on $\partial^{+}\Omega$  is required if $\tau_{1}+2\alpha>\pi.$

Repeat the same argument with $\lambda_{1} = \underline{\gamma}^{+}, \lambda_{2} = \bar{\gamma}^{+}$  and $\gamma_{2} = \bar{\gamma}^{-}$. 
In the same way as above, there exist functions 
\begin{equation}
b_{*}^{+} (x,y) = f(\textbf{w}) + p(\delta) + h_{\beta_{1}}^{-} (x,y) \mathrm{ \ for \ } (x,y) \in \Delta_{\beta_{1}}
\end{equation}
\begin{equation}
b_{*}^{-} (x,y) = f(\textbf{w}) - p(\delta) - h_{\beta_{2}}^{+} (x,y) \mathrm{ \ for \ } (x,y) \in \Delta_{\beta_{2}}
\end{equation} 
such that 
\begin{equation}
\label{R3}
b_{*}^{-}(x,y) < f(x,y)
\end{equation} for $(x,y) \in \Delta_{\beta_{2}} \cap B_{\rho(\delta)}'$ where $B_{\rho(\delta)}' = Y(B_{\rho(\delta)})$. 
Then combining (\ref{R4}) and (\ref{R3}) we get 
\begin{equation}
b_{*}^{-}(x,y) < f(x,y) < b^{+} (x,y)
\end{equation}
for $(x,y) \in \Delta_{\beta_{1}} \cap \Delta_{\beta_{2}} \cap B_{\rho(\delta)}'$.  
As in \cite{NoraKirk2}, it follows that $c(u,v)$  is uniformly continuous on $E.$  
\vspace{3mm}

\noindent {\bf Case (ii):} Case (ii) is simply Case (i) reflected about the $xz$-plane and the proof follows as above.
\vspace{3mm}

\begin{figure}[ht]
\label{FigTWO}
\centering
\includegraphics{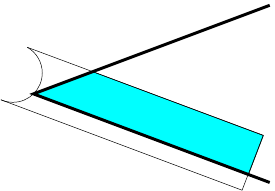}
\caption{The domain of a supersolution in Case (i). \label{FigTwo} }
\end{figure}

\begin{figure}[ht]
\centering
\includegraphics{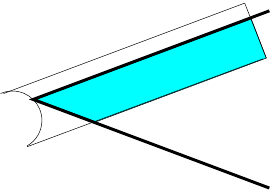}
\caption{The domain of a supersolution in Case (ii). \label{FigThree} }
\end{figure}

\noindent {\bf Case (iii):}  Notice that 
\[ 
0\le \pi-2\alpha < \underline{\gamma}^{+} + \underline{\gamma}^{-} \le \overline{\gamma}^{+}+\overline{\gamma}^{-} < \pi+2\alpha\le\pi 
\]
and so $\underline{\gamma}^{-} = 0$  and $\underline{\gamma}^{+} = 0$  cannot both occur and $\overline{\gamma}^{+}=\pi$  and 
$\overline{\gamma}^{-}=\pi$  cannot both occur.  The result follows from this, using the arguments in Cases 1 \& 2. 
In particular, if $\underline{\gamma}^{-}>0,$  then we obtain a supersolution $ b^{+}$  as in Case (i) (see Figure \ref{FigTwo})
and if $\underline{\gamma}^{-} = 0,$  we obtain a supersolution $ b^{+}$  as in Case (ii) (see Figure \ref{FigThree}); 
if  $\overline{\gamma}^{-}<\pi,$ 
we obtain a subsolution $b_{*}^{-}$  as in Case (i) and if $\overline{\gamma}^{-}=\pi,$  we obtain a subsolution as in Case (ii). 
\vspace{5mm}

\noindent {\bf Now assume (B) holds.} Let $B = \{(x,y) \in\mathbb{R}^{2}: \sqrt{x^{2} + y^{2}} < 1, y \geq 0\}$ and 
let $\bar{B}$ be the closure of $B$ in $\mathbb{R}^{2}$. 
Let $g:\overline{B}\to \overline{E}$  be a conformal or anticonformal map taking $\{(u,0):-1\le u\le 1\}$  onto $\partial E\setminus \sigma$ 
such that the map $X=Y\circ g: B \rightarrow \mathbb{R}^{3}$  has a downward orientation 
(i.e. the normal $X_{u}\times X_{v}$  to $\mathscr{S}_{0}$  gives a downward orientation).  
Writing $X(u,v) = (x(u, v), y(u,v), z(u,v))$  and $K(u, v) = (x(u,v), y(u,v)),$  we have $K\in C^{0}(\bar{B})$  and 
$K(u, 0) = (0,0)$ while $X(u, 0) = (0, 0, z(u, 0))$ for $u \in [-1, 1]$. Then the argument follows from \cite{LS1}  and 
the previous argument here, as explained in \cite{NoraKirk2}.  
  \ \ \ \ $\Box$

\subsection{Proof of Theorem \ref{TWO}:}  
We first note that if $\delta_{1},\delta_{2}\in (-\alpha,\alpha)$  with $\delta_{1}<\delta_{2}$  and 
$Rf(\delta_{1})$  and $Rf(\delta_{2})$  both exist, then it follows from \cite{EL:86}  that $Rf(\theta)$  exists 
for all $\theta\in [\delta_{1},\delta_{2}]$  and  $Rf(\theta)$  is a continuous function of $\theta$  on $[\delta_{1},\delta_{2}]$ 
which behaves as described in {\bf (i)}, {\bf (ii)} or {\bf (iii)} of Theorem \ref{ONE}.  The first part of Theorem \ref{TWO} (i.e. (a)) 
follows from this.

Suppose $\{\epsilon_{j}\}$  is a decreasing sequence with $\lim_{j\to\infty} \epsilon_{j}=0.$   Let $I=(-1,1)$  and set 
\[
\gamma_{j}(s)= \left\{ \begin{array}{ccc} 
\gamma^{+}\left(\epsilon_{j}s\right) & {\rm if} & 0<s<1 \\ 
\gamma^{-}\left(-\epsilon_{j}s\right) & {\rm if} & -1<s<0 \\ 
\end{array}
\right.
\] 
for $j\in \Natural;$   then $\left\{\cos\left(\gamma_{j}\right) : j\in \Natural\right\}$  
is a subset of the unit ball in $L^{\infty}(I)=\left(L^{1}(I)\right)^{*}.$     
By the Banach-Alaoglu theorem, there exist a subsequence $\{\epsilon_{j_{k}}\}$  of $\{\epsilon_{j}\}$   and 
a function $h=h_{\{\epsilon_{j_{k}}\}}\in L^{\infty}(I)$  such that $\cos\left(\gamma_{j_{k}}\right)$  converges weak-star  to $h;$  
that is, for each $m\in L^{1}(I),$    
\[
\lim_{k\to\infty} \int_{-1}^{1} \cos\left(\gamma_{j_{k}}(s)\right)m(s) \ ds =  \int_{-1}^{1} h(s)m(s) \ ds.
\]
Let us define $\gamma^{*}=\gamma^{*}_{\{\epsilon_{j_{k}}\}} =\cos^{-1}\left(h\right)$  (almost everywhere on $(-1,1)$).  
For any $b\in (0,1),$  by choosing $m$  to be the characteristic function of the interval $(0,b)$  we see that  
\[
\int_{0}^{b} h(s) \ ds = \lim_{k\to\infty} \int_{0}^{b} \cos\left(\gamma_{j_{k}}(s)\right) \ ds
= \lim_{k\to\infty} \frac{1}{\epsilon_{j_{k}}} \int_{0}^{b\epsilon_{j_{k}}} \cos\left(\gamma^{+}(t)\right) \ dt 
\]
and, by choosing $m$  to be the characteristic function of the interval $(-b,0),$  we see that
\[
\int_{-b}^{0} h(s) \ ds = \lim_{k\to\infty} \int_{-b}^{0} \cos\left(\gamma_{j_{k}}(s)\right) \ ds
= \lim_{k\to\infty} \frac{1}{\epsilon_{j_{k}}} \int^{b\epsilon_{j_{k}}}_{0} \cos\left(\gamma^{-}(t)\right) \ dt;
\]
hence  
\[
\int_{0}^{b} \cos(\gamma^{*}(s)) \ ds  =  
\lim_{k\to\infty} \frac{1}{\epsilon_{j_{k}}} \int_{0}^{b\epsilon_{j_{k}}} \cos\left(\gamma^{+}(t)\right) \ dt
\ge \liminf_{\epsilon\to 0} \frac{1}{\epsilon}\int_{0}^{\epsilon b} \cos\left(\gamma^{+}(t)\right) dt
\] 
and
\[
\int^{0}_{-b} \cos(\gamma^{*}(s)) \ ds  =  
\lim_{k\to\infty} \frac{1}{\epsilon_{j_{k}}} \int_{0}^{b\epsilon_{j_{k}}} \cos\left(\gamma^{-}(t)\right) \ dt
\ge \liminf_{\epsilon\to 0} \frac{1}{\epsilon}\int_{0}^{\epsilon b} \cos\left(\gamma^{-}(t)\right) dt.
\] 
Thus 
\begin{equation}
\label{first}
\liminf_{\epsilon\to 0} \frac{1}{\epsilon}\int_{0}^{\epsilon b} \cos\left(\gamma^{\pm}(t)\right) dt 
\le \liminf_{j\to \infty} \int_{0}^{b}   \cos\left(\gamma^{\pm}(\epsilon_{j}s)\right) \ ds.   
\end{equation}
Choose a sequence $\{\epsilon_{j}\}$  with $\lim_{j\to\infty} \epsilon_{j}=0$  such that 
\[
\lim_{j\to\infty} \frac{1}{\epsilon_{j}} \int_{0}^{b\epsilon_{j}} \cos\left(\gamma^{+}(t)\right) \ dt 
= \liminf_{\epsilon\to 0} \frac{1}{\epsilon}\int_{0}^{\epsilon b} \cos\left(\gamma^{+}(t)\right) dt;
\] 
as above, there exist a subsequence $\{\epsilon_{j_{k}}\}$  of $\{\epsilon_{j}\}$  and $\gamma_{*}\in L^{\infty}(I)$  
such that $\cos(\gamma_{j_{k}})$  converges weak-star  to $\cos(\gamma_{*}).$   
Then 
\[
\liminf_{\epsilon\to 0} \frac{1}{\epsilon}\int_{0}^{\epsilon b} \cos\left(\gamma^{+}(t)\right) dt 
= \lim_{k\to\infty} \frac{1}{\epsilon_{j_{k}}} \int_{0}^{b\epsilon_{j_{k}}} \cos\left(\gamma^{+}(t)\right) \ dt 
\]
\[
= \lim_{k\to \infty} \int_{0}^{b}   \cos\left(\gamma_{j_{k}}(s)\right) \ ds
= \int_{0}^{b} \cos(\gamma_{*}(s)) \ ds. 
\]

\vspace{5mm}

\noindent {\bf Case 1:}  Suppose case (I) or (DI) of Theorem 1 holds and $\alpha_{2}=\alpha-\beta^{+}.$  
Let us assume there exists  $\lambda\in (\beta^{+},\pi)$  such that 
\begin{equation}
\label{assumption1}
A_{I}^{+}\left(\frac{\sin(\lambda-\beta^{+})}{\sin(\lambda)}\right) + \frac{\sin(\beta^{+})}{\sin(\lambda)} < 1;
\end{equation}
we shall show that this leads to a contradiction.  Set 
\[
b=\frac{\sin(\lambda-\beta^{+})}{\sin(\lambda)}.
\]
Choose a sequence $\{\epsilon_{j}\}$  with $\lim_{j\to\infty} \epsilon_{j}=0$  such that 
\begin{equation}
\label{consequence1}
\lim_{j\to\infty} \frac{1}{\epsilon_{j}} \int_{0}^{b\epsilon_{j}} \cos\left(\gamma^{+}(t)\right) \ dt 
= \liminf_{\epsilon\to 0} \frac{1}{\epsilon}\int_{0}^{\epsilon b} \cos\left(\gamma^{+}(t)\right) dt;
\end{equation}
as before, there exist a subsequence $\{\epsilon_{j_{k}}\}$  of $\{\epsilon_{j}\}$  and $\gamma_{*}\in L^{\infty}(I)$  
such that $\cos(\gamma_{j_{k}})$  converges weak-star  to $\cos(\gamma_{*}).$   
Then 
\[
\liminf_{\epsilon\to 0} \frac{1}{\epsilon}\int_{0}^{\epsilon b} \cos\left(\gamma^{+}(t)\right) dt 
= \lim_{k\to\infty} \frac{1}{\epsilon_{j_{k}}} \int_{0}^{b\epsilon_{j_{k}}} \cos\left(\gamma^{+}(t)\right) \ dt 
\]
\[
= \lim_{k\to \infty} \int_{0}^{b}   \cos\left(\gamma_{j_{k}}(s)\right) \ ds
= \int_{0}^{b} \cos(\gamma_{*}(s)) \ ds. 
\]
Let $\theta_{0} \in (\sigma,\alpha_{2}),$   where $\sigma=\alpha_{1}$  if case (I) holds and $\sigma=\alpha_{R}$  
if case (DI) holds,  and $z_{0}=Rf\left(\theta_{0}\right).$  
Set $\Omega_{k} = \{(x,y) \in \Real^{2} : ( \epsilon_{j_{k}} {\it x}, \epsilon_{j_{k}} {\it y}) \in \Omega \}$ 
and define $f_{k} \in C^{\infty} (\Omega_{k})$ by
\[
f_{k} (x,y) = \frac{1}{\epsilon_{j_{k}}} (f(\epsilon_{j_{k}} x, \epsilon_{j_{k}}y) - z_{0})
\]
for $(x,y) \in \Omega_{k}$.  Let $\gamma_{k}$ be defined on $\partial \Omega_{k} \backslash \{{\cal O}\}$ by 
\[
\gamma_{k}(x,y)=\gamma(\epsilon_{j_{k}}{\it x}, \epsilon_{j_{k}} y)
\] 
and let $\nu_{k} = \nu_{k} (x,y)$ denote the outward unit normal to $\partial \Omega_{k}$.  Then
$f_{k}$ satisfies the capillary problem
\[
\begin{array}{lll}
N f_{k} (x,y) & = & \epsilon_{j_{k}} \kappa f(\epsilon_{j_{k}} x, \, \epsilon_{j_{k}} y) + \epsilon_{j_{k}} \lambda \, ,
\hspace*{.65in} (x,y) \in \Omega_{k} \\[2mm] T f_{k} \cdot \nu_{k} & = & \cos (\gamma_{k}) \hspace*{1.7in}
\mbox{on} \,\, \partial \Omega_{k} \backslash \{{\cal O}\} \, .
\end{array}
\]
Since $Rf(\theta)<z_{0}$  if $\sigma<\theta<\theta_{0}$  and $Rf(\theta)>z_{0}$  if $\theta_{0}<\theta<\alpha,$  
we see (e.g. \cite{CFC,NCFC,Sim:80}; also see \cite{Tam:84,Tam:86b}) 
that $\{f_{k}\}$   converges locally to the generalized solution $f_{\infty}$  (in the sense of Miranda \cite{Mir:77} and Guisti \cite{Gui:80,Gui:84}) 
of the functional
\[
{\cal F}_{\infty} (g) = \int \int\limits_{\Omega_{\infty}} \sqrt{1+ \mid Dg \mid^{2}} \; dx  
- \int_{\partial\Omega_{\infty}} \cos\left(\gamma_{*}(s)\right) g \; ds, 
\]
where 
\[
f_{\infty}(r\cos(\theta),r\sin(\theta))= \left\{ \begin{array}{ccc} 
-\infty & {\rm if} & -\alpha<\theta<\theta_{0} \\ 
\infty & {\rm if} & \theta_{0}<\theta<\alpha \\ \end{array}
\right.
\] 
if case (I) holds or case (DI) holds and $z_{0}> Rf(\theta)$  for all $\theta\in (-\alpha,\alpha_{L})$  and 
\[
f_{\infty}(r\cos(\theta),r\sin(\theta))= \left\{ \begin{array}{ccc} 
\infty & {\rm if} & -\alpha<\theta<\theta_{h} \\
-\infty & {\rm if} & \theta_{h}<\theta<\theta_{0} \\ 
\infty & {\rm if} & \theta_{0}<\theta<\alpha \\ \end{array}
\right.
\] 
with $Rf\left(\theta_{h}\right)=z_{0}$  and $\theta_{h}<\alpha_{L}$  otherwise.  

\begin{figure}[ht]
\label{FigONE}
\centering
\includegraphics{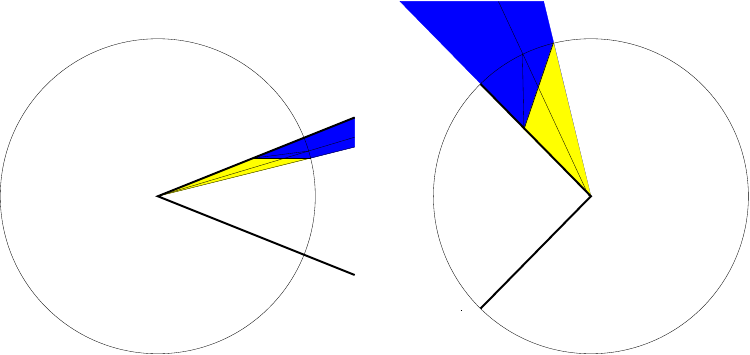}
\caption{The yellow region represents $\Sigma_{\theta_{0}}.$ \label{FigOne} }
\end{figure}

Let us now define the sets
$$
{\cal P}=\{(x,y)\in\Omega_{\infty}:f_{\infty}(x,y)=\infty\}\ \  
{\rm and} \ \ 
{\cal N}=\{(x,y)\in\Omega_{\infty}:f_{\infty}(x,y)=-\infty\}.
$$
These sets have a special structure which follows from the fact that ${\cal P}$ minimizes the functional
\[
\Phi (A)  =  \int \int\limits_{\Omega_{\infty}} \mid D \chi_{A}\mid - \int_{\partial\Omega_{\infty}} \cos(\gamma_{*})\chi_{A} d H^{1}  
\]
and ${\cal N}$ minimizes the functional
\[
\Psi (A)  =  \int \int\limits_{\Omega_{\infty}} \mid D \chi_{A} \mid 
+ \int_{\partial\Omega_{\infty}} \cos(\gamma_{*})\chi_{A} d H^{1}  
\]
in the appropriate sense (e.g. \cite{Gui:80,LS1,Mir:77}).  
Let $\Sigma_{\theta_{0}}$  denote the (open) triangular region whose boundary is the triangle with vertices $(0,0),$  
$B=(b\cos(\alpha),b\sin(\alpha))$  and $C=\left(\cos(\theta_{0}),\sin(\theta_{0})\right)$   and set 
$A={\cal P}\setminus \Sigma_{\theta_{0}}$  (see Figure \ref{FigOne}).  
Simple trigonometric computations with $R>2$ show that
\begin{equation}
\label{trig1}
\Phi\left(B({\cal O},R)\cap {\cal P}\right) -  \Phi\left(B({\cal O},R)\cap {\cal P}\setminus \Sigma_{\theta_{0}}\right) 
= \left(1 - A_{I}^{+}(b)\right) - \left(\frac{\sin(\alpha-\theta_{0})}{\sin(\omega)}\right), 
\end{equation}
where $\pi-\omega$  is the angle $\angle {\cal O}BC.$  
This holds for all $\theta_{0}<\alpha_{2}=\alpha-\beta^{+};$  taking the limit as $\theta_{0}\uparrow \alpha-\beta^{+}$  
and noticing that $\omega\to\lambda$  as $\theta_{0}\uparrow \alpha-\beta^{+},$  we see that 
\[
\Phi\left(B({\cal O},R)\cap {\cal P}\right) -  \Phi\left(B({\cal O},R)\cap {\cal P}\setminus \Sigma_{\alpha_{2}}\right) 
= \left(1 - A_{I}^{+}(b)\right) - \left(\frac{\sin(\beta^{+})}{\sin(\lambda)}\right) > 0
\]
or 
\[
\Phi\left(B({\cal O},R)\cap {\cal P}\right) > \Phi\left(B({\cal O},R)\cap {\cal P}\setminus \Sigma_{\alpha_{2}}\right);
\]
this contradicts the fact that ${\cal P}$  (locally) minimizes $\Phi.$   Therefore (\ref{assumption1}) is false.  
This completes Case 1. 
\vspace{5mm}

\noindent  {\bf Case 2.} Suppose case (D) or (DI) of Theorem 1 holds and $\alpha_{1}=-\alpha+\beta^{-}.$  
Let us assume there exists  $\lambda\in (\beta^{-},\pi)$  such that 
\begin{equation}
\label{assumption2}
A_{I}^{-}\left(\frac{\sin(\lambda-\beta^{-})}{\sin(\lambda)}\right) + \frac{\sin(\beta^{-})}{\sin(\lambda)} < 1.
\end{equation}  
Using a similar argument to that in Case 1, we reach a contradiction. 
\vspace{5mm}

\noindent  {\bf Case 3.} Suppose case (I) or (ID) of Theorem 1 holds and $\alpha_{1}=-\alpha+\beta^{-}.$  
Let us assume there exists  $\lambda\in (\beta^{-},\pi)$  such that 
\begin{equation}
\label{assumption3}
1+A_{S}^{-}\left(\frac{\sin(\lambda-\beta^{-})}{\sin(\lambda)}\right)  > \frac{\sin(\beta^{-})}{\sin(\lambda)}.
\end{equation}
Set 
\[
b=\frac{\sin(\lambda-\beta^{-})}{\sin(\lambda)}.  
\]
Arguing as in Case 1., we see that the set  ${\cal N}=\{(x,y)\in\Omega_{\infty}:f_{\infty}(x,y)=-\infty\}$  minimizes the functional
\[
\Psi (A)  =  \int \int\limits_{\Omega_{\infty}} \mid D \chi_{A} \mid 
+ \int_{\partial\Omega_{\infty}} \cos(\gamma_{*})\chi_{A} d H^{1}  
\]
in the appropriate sense (e.g. \cite{Gui:80,LS1,Mir:77}).  
Let $\Sigma_{\theta_{0}}$  denote the (open) triangular region whose boundary is the triangle with vertices $(0,0),$  
$B=(b\cos(-\alpha),b\sin(-\alpha))$  and $C=\left(\cos(\theta_{0}),\sin(\theta_{0})\right)$   and set 
$A={\cal N}\setminus \Sigma_{\theta_{0}}$  (see Figure \label{FigTHREE}).  
Simple trigonometric computations with $R>2$ show that\[
\Psi\left(B({\cal O},R)\cap {\cal N}\right) -  \Psi\left(B({\cal O},R)\cap {\cal N}\setminus \Sigma_{\theta_{0}}\right) 
= \left(1 + A_{S}^{-}(b)\right) - \left(\frac{\sin(\alpha+\theta_{0})}{\sin(\omega)}\right), 
\]
where $\omega$  is the angle $\angle {\cal O}BC.$  
This holds for all $\theta_{0}>\alpha_{1}=-\alpha+\beta^{-};$  taking the limit as $\theta_{0}\downarrow -\alpha+\beta^{-}$  
and noticing that $\omega\to\lambda$  as $\theta_{0}\downarrow -\alpha+\beta^{-},$  we see that 
\[
\Psi\left(B({\cal O},R)\cap {\cal N}\right) -  \Psi\left(B({\cal O},R)\cap {\cal N}\setminus \Sigma_{\alpha_{1}}\right) 
= \left(1 + A_{S}^{-}(b)\right) - \left(\frac{\sin(\beta^{-})}{\sin(\lambda)}\right) > 0
\]
or 
\[
\Psi\left(B({\cal O},R)\cap {\cal N}\right) > \Psi\left(B({\cal O},R)\cap {\cal N}\setminus \Sigma_{\alpha_{1}}\right);
\]
this contradicts the fact that ${\cal N}$  (locally) minimizes $\Psi.$   Therefore (\ref{assumption3}) is false.  
This completes Case 3. 
\vspace{5mm}

\noindent  {\bf Case 4.} Suppose case (D) or (ID) of Theorem 1 holds and $\alpha_{2}=\alpha-\beta^{+}.$  
Let us assume there exists  $\lambda\in (\beta^{+},\pi)$  such that 
\begin{equation}
\label{assumption4}
1+A_{S}^{+}\left(\frac{\sin(\lambda-\beta^{+})}{\sin(\lambda)}\right)  > \frac{\sin(\beta^{+})}{\sin(\lambda)}.
\end{equation}
Using a similar argument to that in Case 3, we reach a contradiction. 
The proof of Theorem \ref{TWO} is then complete. \ \ \ \ $\Box$

\section{Corollaries and Examples} 

\begin{cor}  
\label{COR1}
Suppose $m\in [-1,1];$   set $\sigma=\cos^{-1}(m)\in [0,\pi].$ 
\begin{itemize}
\item[(a)]  If $A^{+}_{I}(b)\le mb$  and case (I) or (DI) holds, then $\beta^{+}\ge \sigma.$     
\item[(b)]  If $A^{-}_{I}(b)\le mb$  and case (D) or (DI) holds, then $\beta^{-}\ge \sigma.$    
\item[(c)]  If $A^{-}_{S}(b)\ge mb$  and case (I) or (ID) holds, then $\beta^{-}\ge \pi-\sigma.$    
\item[(d)]  If $A^{+}_{S}(b)\ge mb$  and case (D) or (ID) holds, then $\beta^{+}\ge \pi-\sigma.$    
\end{itemize}
\end{cor}

\noindent {\bf Proof:}  (a)  Suppose case (I) or (DI) of Theorem \ref{ONE} holds, $\sigma\in [0,\pi],$  $\cos(\sigma)=m,$  and $\beta^{+}<\sigma.$  
By Theorem \ref{TWO} (a), we know that 
\[
\sin(\sigma)\left(\frac{\sin(\lambda-\beta^{+})}{\sin(\lambda)}\right) + \frac{\sin(\beta^{+})}{\sin(\lambda)}\ge 
A_{I}^{+}\left(\frac{\sin(\lambda-\beta^{+})}{\sin(\lambda)}\right) + \frac{\sin(\beta^{+})}{\sin(\lambda)}\ge 1
\]
or 
\[
\frac{\cos(\sigma)\sin(\lambda-\beta^{+}) + \sin(\beta^{+})}{\sin(\lambda)}\ge 1
\]
for all $\lambda\in (\beta^{+},\pi).$   Since $\sigma>\beta^{+},$  we may set $\lambda=\sigma$  and obtain 
\[
\cos(\sigma-\beta^{+})=\frac{\cos(\sigma)\sin(\sigma-\beta^{+}) + \sin(\beta^{+})}{\sin(\sigma)}\ge 1, 
\]
which is a contradiction since $\sigma-\beta^{+}\neq 0.$  Thus $\beta^{+}\ge \sigma.$

\noindent (b) This is essentially the same as (a). 

\noindent (c) Suppose case (I) or (ID) of Theorem \ref{ONE} holds, $\sigma\in [0,\pi],$  $\cos(\sigma)=m,$  and $\beta^{-}<\pi-\sigma.$  
By Theorem \ref{TWO} (c), we know that 
\[
1+\sin(\sigma)\left(\frac{\sin(\lambda-\beta^{-})}{\sin(\lambda)}\right)  \le 
1+A_{S}^{-}\left(\frac{\sin(\lambda-\beta^{-})}{\sin(\lambda)}\right)  \le \frac{\sin(\beta^{-})}{\sin(\lambda)}
\]
or
\[
\frac{\sin(\lambda)+\cos(\sigma)\sin(\lambda-\beta^{-}) - \sin(\beta^{-})}{\sin(\lambda)}\le 0
\]
for all $\lambda\in (\beta^{-},\pi).$   Since $\beta^{-}<\pi-\sigma,$  we may set $\lambda=\pi-\sigma$  and obtain 
\[
1+\cos\left(\sigma+\beta^{-}\right) = \frac{\sin(\sigma)+\cos(\sigma)\sin(\sigma+\beta^{-}) - \sin(\beta^{-})}{\sin(\sigma)} \le 0, 
\]
which is a contradiction since $\sigma+\beta^{-}<\pi.$  Thus $\beta^{-}\ge \pi-\sigma.$

\noindent (d)  This is essentially the same as (c).  \hfill $\Box$
\vspace{5mm}

\begin{example}
\label{ex1}
Let $\alpha\in (0,\pi]$  and $\gamma^{\pm}_{1},\gamma^{\pm}_{2}\in [0,\pi]$   with $\gamma^{+}_{1}\le \gamma^{+}_{2}$  
and $\gamma^{-}_{1}\le \gamma^{-}_{2}.$
Set 
\[
\Omega=\{(r\cos(\theta),r\sin(\theta)) : 0<r<1, -\alpha<\theta<\alpha \}.
\]
For each $n\in\Natural,$  let $A_{n}=\left(2^{-n^2}, 2^{-n(n-1)}\right]$  and  $B_{n}= \left(2^{-n(n+1)},2^{-n^2}\right].$
Define 
\[
\gamma(s)= \sum_{n=1}^{\infty} \left(\gamma^{+}_{1}I_{A_{n}}(s)+ \gamma^{+}_{2}I_{B_{n}}(s)
+\gamma^{-}_{1}I_{A_{n}}(-s)+ \gamma^{-}_{2}I_{B_{n}}(-s)\right), 
\]
so that $\gamma$  is defined on $\partial\Omega\cap B\left({\cal O},1\right)$  by  
\[
\gamma(r\cos(\theta),r\sin(\theta))= \left\{ \begin{array}{ccc} 
\gamma^{+}_{1} & {\rm if} & \theta=\alpha,\ 2^{-n^2}<r\le 2^{-n(n-1)} \ \ {\rm for \ some} \ n\in\Natural\\ 
\gamma^{+}_{2} & {\rm if} & \theta=\alpha,\ 2^{-n(n+1)}<r\le 2^{-n^2} \ \ {\rm for \ some} \ n\in\Natural\\
\gamma^{-}_{1} & {\rm if} & \theta=-\alpha,\ 2^{-n^2}<r\le 2^{-n(n-1)} \ \ {\rm for \ some} \ n\in\Natural\\ 
\gamma^{-}_{2} & {\rm if} & \theta=-\alpha,\  2^{-n(n+1)}<r\le 2^{-n^2} \ \ {\rm for \ some} \ n\in\Natural.\\
\end{array}
\right.
\]
Set 
\[
c_{j}=\left\{ \begin{array}{ccc}
2^{-\frac{j}{2}\left(\frac{j}{2}+1\right)}  & {\rm if} & j \ {\rm is \ even}\\ 
2^{-\left(\frac{j+1}{2}\right)^2}  & {\rm if} & j \ {\rm is \ odd}.\\ 
\end{array}
\right.
\]
Let $b\in (0,1)$  be fixed for now.  Set $\epsilon_{j}=\frac{c_{2j}}{b}$  ($j\in\Natural$); notice that 
$\frac{c_{2j+1}}{c_{2j}}=2^{-(j+1)}.$   Then 
\[
b\cos(\gamma_{1}^{\pm})\ge A_{S}^{\pm}(b) \ge 
\lim_{j\to\infty} \frac{1}{\epsilon_{j}} \int_{0}^{\epsilon_{j}b} \cos(\gamma^{\pm}(t)) \ dt 
= \lim_{j\to\infty} b \int_{0}^{1} \cos(\gamma^{\pm}_{j}(sb)) \ ds
\]
\[
= \lim_{j\to\infty} b \int_{0}^{1} \cos(\gamma^{\pm}(c_{2j}s)) \ ds 
\]
\[
= \lim_{j\to\infty} b \left(\int_{\frac{c_{2j+1}}{c_{2j}}}^{1}  \cos(\gamma^{\pm}(c_{2j}s)) \ ds 
+ \int_{0}^{\frac{c_{2j+1}}{c_{2j}}} \cos(\gamma^{\pm}(c_{2j}s)) \ ds \right)
\]
\[
= \lim_{j\to\infty} b \left(\cos(\gamma^{\pm}_{1}) \left(1-2^{-(j+1)}\right)
+ \int_{0}^{2^{-(j+1)}} \cos(\gamma^{\pm}(c_{2j}s)) \ ds \right) 
= b\cos(\gamma^{\pm}_{1}).
\]
Using a similar argument for $A^{\pm}_{I}(b)$  with $\epsilon_{j}=\frac{c_{2j+1}}{b},$  $j\in\Natural,$ 
we see that  
\begin{equation}
\label{bound1}
A^{\pm}_{I}(b) = b\cos(\gamma^{\pm}_{2})  \ \ \ \ {\rm and}  \ \ \ \ A^{\pm}_{S}(b) = b\cos(\gamma^{\pm}_{1}).
\end{equation} 
\end{example}

\vspace{5mm}

\begin{example} 
\label{ex2}
Let $\alpha\in (0,\pi]$  and $\gamma^{\pm}_{1},\gamma^{\pm}_{2}\in [0,\pi]$   with $\gamma^{+}_{1}\le \gamma^{+}_{2}$  
and $\gamma^{-}_{1}\le \gamma^{-}_{2}.$
Set 
\[
\Omega=\{(r\cos(\theta),r\sin(\theta)) : 0<r<1, -\alpha<\theta<\alpha \}.
\]
For each $n\in\Natural,$  let $A_{n}=\left(\frac{2}{4^{n}},\frac{4}{4^{n}}\right),$   $B_{n}=\left(\frac{1}{4^{n}},\frac{2}{4^{n}}\right),$  
and $C_{n}=\{\frac{4}{4^{n}}\}.$ 
Define 
\[
\gamma(s)= \sum_{n=1}^{\infty} \left(\gamma^{+}_{1}I_{A_{n}}(s)+ \gamma^{+}_{2}I_{B_{n}}(s)+\pi I_{C_{n}}(s)
+\gamma^{-}_{1}I_{A_{n}}(-s)+ \gamma^{-}_{2}I_{B_{n}}(-s)+\pi I_{C_{n}}(-s) \right), 
\]
so that $\gamma$  is defined on $\partial\Omega\cap B\left({\cal O},1\right)$  by  
\[
\gamma(r\cos(\theta),r\sin(\theta))= \left\{ \begin{array}{ccc} 
\gamma^{+}_{1} & {\rm if} & \theta=\alpha,\ \frac{2}{4^{n}}<r<\frac{4}{4^{n}} \ \ {\rm for \ some} \ n\in\Natural\\ 
\gamma^{+}_{2} & {\rm if} & \theta=\alpha,\ \frac{1}{4^{n}}<r<\frac{2}{4^{n}} \ \ {\rm for \ some} \ n\in\Natural\\
\pi & {\rm if} & \theta=\alpha,\ r=\frac{4}{4^{n}} \ \ {\rm for \ some} \ n\in\Natural\\ 
0 & {\rm if} & \theta=\alpha,\ r=\frac{2}{4^{n}} \ \ {\rm for \ some} \ n\in\Natural\\ 
\gamma^{-}_{1} & {\rm if} & \theta=-\alpha,\ \frac{2}{4^{n}}<r<\frac{4}{4^{n}} \ \ {\rm for \ some} \ n\in\Natural\\ 
\gamma^{-}_{2} & {\rm if} & \theta=-\alpha,\ \frac{1}{4^{n}}<r<\frac{2}{4^{n}} \ \ {\rm for \ some} \ n\in\Natural\\
\pi & {\rm if} & \theta=-\alpha,\ r=\frac{4}{4^{n}} \ \ {\rm for \ some} \ n\in\Natural\\ 
0 & {\rm if} & \theta=-\alpha,\ r=\frac{2}{4^{n}} \ \ {\rm for \ some} \ n\in\Natural.\\ 
\end{array}
\right.
\]
Then 
\[
\liminf_{r\to 0} \gamma\left(r\cos(\pm \alpha),r\sin(\pm \alpha)\right)=0,   \ \ 
\limsup_{r\to 0} \gamma\left(r\cos(\pm \alpha),r\sin(\pm \alpha)\right)=\pi, 
\]
\[
\essliminf_{r\to 0} \gamma\left(r\cos(\pm \alpha),r\sin(\pm \alpha)\right)=\gamma^{\pm}_{1}, \   
\esslimsup_{r\to 0} \gamma\left(r\cos(\pm \alpha),r\sin(\pm \alpha)\right)=\gamma^{\pm}_{2}.
\]
Thus $A^{\pm}_{I}(b)\ge b\cos\left(\gamma^{\pm}_{2}\right)$  and $A^{\pm}_{S}(b)\le b\cos\left(\gamma^{\pm}_{1}\right).$  

Let $b\in (0,1)$  be fixed for now.  If we set $\epsilon_{j}=\frac{1}{b4^{j}}$  ($j\in\Natural$), 
then $\gamma_{j}\left(s\right)=\gamma\left(\frac{s}{b}\right)$  and so  
\[
A^{+}_{S}(b) \ge \lim_{j\to\infty} \int_{0}^{b} \cos\left(\gamma_{j}(s)\right) \ ds = b\int_{0}^{1} \cos\left(\gamma(s)\right) \ ds 
=b\left(\frac{2}{3}\cos(\gamma^{+}_{1})+\frac{1}{3}\cos(\gamma^{+}_{2})\right)
\]
and, if  we set $\epsilon_{j}=\frac{2}{b4^{j}}$  ($j\in\Natural$), then $\gamma_{j}\left(s\right)=\gamma\left(\frac{s}{2b}\right)$  and so  
\[
A^{+}_{I}(b) \le \lim_{j\to\infty} \int_{0}^{b} \cos\left(\gamma_{j}(s)\right) \ ds = 2b\int_{0}^{\frac{1}{2}} \cos\left(\gamma(s)\right) \ ds 
=b\left(\frac{1}{3}\cos(\gamma^{+}_{1})+\frac{2}{3}\cos(\gamma^{+}_{2})\right);
\]
similar estimates hold on $\partial^{-}\Omega.$  
Now suppose $\left(\eta_{j}\right)$  is any decreasing sequence in $(0,1)$  converging to zero.  
For each $j\in\Natural,$  there exists a $k \in\Natural$  such that $\frac{1}{4} \le 4^{k-1}\eta_{j}b<1$  and, since 
$\gamma$  is piecewise constant, a direct calculation shows that 
$\int_{0}^{b}   \cos\left(\gamma^{\pm}(\eta_{j}s)\right) \ ds = b \int_{0}^{1} \cos\left(\gamma^{\pm}(\eta_{j}bs)\right) \ ds$  equals  
\[
\left\{ \begin{array}{ccc} 
b\left( \frac{1}{6\cdot 4^{k-1}\eta_{j}b} \cos\gamma_{1}^{\pm} + \left(1-\frac{1}{6\cdot 4^{k-1}\eta_{j}b}\right) \cos\gamma_{2}^{\pm} \right)
& {\rm if} &  \frac{1}{4} \le 4^{k-1}\eta_{j}b<\frac{1}{2}      \\ 
b\left(\left(1-\frac{1}{3\cdot 4^{k-1}\eta_{j}b}\right) \cos\gamma_{1}^{\pm} + \frac{1}{3\cdot 4^{k-1}\eta_{j}b} \cos\gamma_{2}^{\pm} \right)
& {\rm if} & \frac{1}{2} \le 4^{k-1}\eta_{j}b<1.\\ 
\end{array}
\right.
\]
The minumum occurs when $4^{k-1}\eta_{j}b=\frac{1}{2}$  and the minimum of  $\int_{0}^{b}   \cos\left(\gamma^{\pm}(\eta_{j}s)\right) ds$  is 
$b\left(\frac{1}{3}\cos(\gamma^{\pm}_{1})+\frac{2}{3}\cos(\gamma^{\pm}_{2})\right).$   
The maximum occurs when $4^{k-1}\eta_{j}b=\frac{1}{4}$  and the maximum of  $\int_{0}^{b}   \cos\left(\gamma^{\pm}(\eta_{j}s)\right) ds$  is 
$b\left(\frac{2}{3}\cos(\gamma^{\pm}_{1})+\frac{1}{3}\cos(\gamma^{\pm}_{2})\right).$   
Thus 
\begin{equation}
\label{bound2}
A^{\pm}_{I}(b) = b\left(\frac{1}{3}\cos(\gamma^{\pm}_{1})+\frac{2}{3}\cos(\gamma^{\pm}_{2})\right) 
\end{equation} 
and 
\begin{equation}
\label{bound3}
A^{\pm}_{S}(b) = b\left(\frac{2}{3}\cos(\gamma^{\pm}_{1})+\frac{1}{3}\cos(\gamma^{\pm}_{2})\right).
\end{equation} 
\end{example}
\vspace{5mm}

\noindent 
In these Examples, we have the same essential limits inferior and superior at ${\cal O}$  and yet $A^{\pm}_{I}$  and $A^{\pm}_{S}$  
behave differently.  In Example \ref{ex1}, we have the ``extreme values'' (\ref{bound1}); the ``effective'' contact angles in 
(a) and (b) of Corollary \ref{COR1} are $\gamma^{\pm}_{2}$  and in (c) and (d) of Corollary \ref{COR1} are $\gamma^{\pm}_{1}.$ 
On the other hand, in  Example \ref{ex2}, we have the ``intermediate values'' (\ref{bound2}) and (\ref{bound3}). 
For Example \ref{ex2}, the ``effective'' contact angles in (a) and (b) of Corollary \ref{COR1} are $\sigma_{2}^{\pm}$   
and in (c) and (d) of Corollary \ref{COR1} are $\sigma_{1}^{\pm},$  where $\sigma_{1}^{\pm},\sigma_{2}^{\pm}\in [0,\pi]$  
satisfy  
\[
\cos \sigma_{1}^{\pm}= \frac{2}{3}\cos\gamma^{\pm}_{1}+\frac{1}{3}\cos\gamma^{\pm}_{2} \ \ \ \ {\rm and} \ \ \ \ 
\cos \sigma_{2}^{\pm}= \frac{1}{3}\cos\gamma^{\pm}_{1}+\frac{2}{3}\cos\gamma^{\pm}_{2}.
\]
If $f$  is a bounded solution of (\ref{eq:abasic}) satisfying (\ref{eq:bbasic}) on $\partial^{\pm}\Omega^{*} \setminus \{{\cal O}\}$ 
which is discontinuous at ${\cal O}$  and  $Rf(\theta)$ exists for all $\theta \in (-\alpha,\alpha),$  then 
bounds on the sizes $\beta^{+}$  and $\beta^{-}$  of side fans can be computed using Corollary \ref{COR1};  
the lower bounds on the sizes of these side fans differ between these two examples. 

\section{Comments and Extensions}  The last section of \cite{LS1} dealt with extensions of (\ref{eq:abasic}) to equations of 
prescribed mean curvature.  Consider the prescribed mean curvature contact angle problem 
\begin{eqnarray}
\label{PMC}
Nf & = & 2H(\cdot,f) \mbox{  \ in \ } \Omega    \\
Tf \cdot {\bf \nu} & = & \cos \gamma  \ \ \ \mbox{a.e.  on \ } \partial \Omega.
\label{PMC:contact}
\end{eqnarray}
Suppose $f\in C^{2}(\Omega)$  satisfies (\ref{PMC}) and (\ref{PMC:contact}) and also suppose the following conditions hold:  
\begin{itemize}
\item[(i)] $\sup_{x\in\Omega} |f(x)| < \infty \ \ \ \  {\rm and} \ \ \ \  \sup_{x\in\Omega} |H(x,f(x))| < \infty.$  
\item[(ii)] $H(x,y,t)$  is weakly increasing in $t$  for each $(x,y)\in\Omega.$
\end{itemize}
Using \cite{EchartLancaster1}, we see that Theorems \ref{ONE} and \ref{TWO} continue to hold for solutions $f$  as above; 
the argument is the same as that in \cite{LS1}.

One might ask if the case considered in Theorem \ref{TWO} is of ``physical interest.''  Is it possible for the contact angle to fail to have a 
limit at the corner ${\cal O}$?  In a sense this is a silly question since, at a small enough scale, the macroscopic description of a capillary surface
becomes meaningless.  On the other hand, one sometimes uses devices (e.g. homogenization) to obtain useful macroscopic information from knowledge of 
``small scale'' properties.  An experiment which might be of some interest would be to form a vertical wedge consisting of two planes of glass which have 
been coated in increasing narrow vertical strips with a non-wetting substance (e.g. paraffin) as the edge at which the two planes meet is approached; 
this would approximate the situation considered in Theorem \ref{TWO} and one wonders if there is a ``effective'' contact angle at the corner which is 
larger than that for glass and smaller than that for paraffin.

\end{document}